\documentclass{rsauthor}%%%%where rsproca is the template name
\usepackage{graphicx,psfrag,overpic}
\usepackage{amssymb,mathrsfs,amscd}
\usepackage[numbers]{natbib}

\jname{Phil. Trans. R. Soc. A. }

\newtheorem{problemL}{\bf Problem}[section]
\newcommand{\x}{{\bf x}}
\newcommand{\R}{\mathbb{R}}

\newcommand{\vv}{\mathbf{v}}
\newcommand{\aaa}{\mathbf{a}}
\newcommand{\bb}{\mathbf{b}}

\newcommand \alp{\alpha}

\newcommand \vphi{\varphi}

\newcommand \gam{\gamma}

\newcommand \der{\partial}

\newcommand{\PtUpL}{{P_1}}
\newcommand{\PtLwL}{{P_2}}
\newcommand{\PtLwR}{{P_3}}
\newcommand{\PtUpR}{{P_4}}

\newcommand \Wedge{\Gamma_{\rm wedge}}
\newcommand \Om{\Omega}

\newcommand \irho{\rho_{\infty}}
\newcommand \ivphi{\varphi_{\infty}}

\newcommand \iu{u_{\infty}}

\newcommand{\nnu}{{\boldsymbol\nu}}

\newcommand \ol{\overline}

\newcommand \leftbottom{P_4}
\newcommand \rightbottom{P_3}

\newcommand \shock{\Gamma_{\rm shock}}

\newcommand{\mr}{\mathbb{R}}

\newcommand{\deq}{:=}%{\overset{\underset{\mathrm{def}}{}}{=}}

\newcommand{\pSi}{\varphi}

\renewcommand{\a}{\alpha}

\renewcommand{\r}{\rho}

\renewcommand{\t}{\theta}
\newcommand{\e}{\epsilon}

\renewcommand{\O}{\Omega}

\renewcommand{\L}{\Lambda}

\newcommand{\PtIncW}{{P_0}}
\newcommand{\Gwt}{\Gamma_{\text{\rm wedge}}^2}

\newcommand{\Gsh}{\Gamma_{\text{\rm shock}}}
\newcommand \leftshockop{S_{0}}
\newcommand \leftsonicop{\Gamma_{\rm sonic}^{0}}

\newcommand \oOmop{\Om_0}

\newcommand \lefttop{P_1}
\newcommand \righttop{P_2}

\newcommand \rightshockop{S_{1}}
\newcommand \rightsonicop{\Gamma_{\rm sonic}^{1}}

\newcommand \nOmop{\Om_1}

\begin{document}

\title{Free Boundary Problems in Shock Reflection/Diffraction and Related Transonic Flow Problems}

\author{Gui-Qiang G. Chen$^{1}$ and Mikhail Feldman$^{2}$ }

\address{$^{1}$ Mathematical Institute, University of
Oxford,
Oxford, OX2 6GG, UK.
E-mail: chengq@maths.ox.ac.uk\\
$^{2}$ Department of Mathematics, University of Wisconsin,
WI 53706, USA.
E-mail: feldman@math.wisc.edu}

%%%% Subject entries to be placed here %%%%
%\subject{35-02, 35R35, 35R37, 35M12, 35L65, 35L67, 35B30, 35B65, 35J70, 35D30, 76H05, 76L05, 76N10, 35Q35, 35B35, 35B40, 35J25}


\corres{Gui-Qiang G. Chen\\
\email{chengq@maths.ox.ac.uk}}

\abstract
{\noindent Shock waves are steep wave fronts that are fundamental in nature, especially in high-speed fluid flows.
When a shock hits an obstacle, or a flying body meets a shock, shock reflection/diffraction phenomena
occur. In this paper, we show how several longstanding
shock reflection/diffraction
problems can be formulated as free boundary problems,
discuss some recent progress
in developing
mathematical
ideas, approaches, and techniques for solving these
problems, and present some further open problems in this direction.
In particular, these shock problems include
von Neumann's problem for shock reflection-diffraction by two-dimensional wedges with concave corner,
Lighthill's problem for shock diffraction by two-dimensional wedges with convex corner,
and Prandtl-Meyer's problem for supersonic flow
impinging
onto solid wedges,
which are also fundamental in the mathematical theory of multidimensional
conservation laws.
}

%%%%%%%%%%%%%%%%%%%%%%%%%%%
\maketitle
%%%%%%%%%% Insert the texts which can accomdate on firstpage in the tag "fmtext" %%%%%

\section{Introduction}

Shock waves are steep fronts that propagate in the compressible fluids in which convection dominates diffusion.
They are fundamental in nature, especially in high-speed fluid flows.
Examples include transonic and/or supersonic shocks formed by supersonic flows impinging onto solid wedges,
transonic shocks around supersonic or near sonic flying bodies,
bow shocks created by solar winds in space, blast waves by explosions,
and other shocks by various natural processes.
When a shock hits an obstacle, or a flying body meets a shock,
shock reflection/diffraction phenomena occur.

Many of such shock reflection/diffraction problems
can be formulated as free boundary
problems involving nonlinear partial differential equations (PDEs) of mixed
elliptic-hyperbolic type.
The understanding of these shock reflection/diffraction phenomena
requires a complete mathematical solution of the corresponding free
boundary problems for nonlinear mixed PDEs.
In this paper, we show how several longstanding, fundamental multidimensional shock problems
can be formulated as free boundary problems,
discuss some recent progress in developing
mathematical
ideas, approaches, and techniques for solving these
problems, and present some further open problems in this direction.
In particular, these shock problems include
von Neumann's problem for shock reflection-diffraction
by two-dimensional wedges with concave corner,
Lighthill's problem for shock diffraction by two-dimensional wedges
with convex corner, and
Prandtl-Meyer's problem for supersonic flow impinging onto solid wedges.
These problems are not only longstanding open problems in fluid mechanics,
but also fundamental in the mathematical theory of multidimensional
conservation laws:
These shock reflection/diffration configurations are the core configurations
in the structure of global entropy solutions of the two-dimensional
Riemann problem for hyperbolic conservation laws,
while the Riemann solutions are building blocks and local structure of general
solutions and determine global attractors and asymptotic states of entropy
solutions, as time tends to infinity, for multidimensional hyperbolic systems
of conservation laws; see \cite{Chen-Feldman15,Glimm-Majda,Serre07,Zheng01}
and the references cited therein. In this sense,
we have to understand the shock reflection/diffraction phenomena
in order to understand fully
global entropy solutions to multidimensional hyperbolic systems of conservation
laws.

We first focus on these problems for the Euler equations for potential flow.
The unsteady potential flow is governed by the conservation law
of mass and Bernoulli's law:
\begin{align}
\label{1-a}
&\partial_t\rho+\nabla_{\bf x}\cdot(\rho \nabla_{\bf x}\Phi)=0,\\
\label{1-b}
&\partial_t\Phi+\frac 12|\nabla_{\bf x}\Phi|^2+h(\rho)=B
\end{align}
for the density $\rho$ and the velocity potential $\Phi$, where the Bernoulli constant
$B$ is determined by the incoming flow and/or boundary conditions,
and $h(\rho)$ satisfies the relation
\begin{equation*}
h'(\rho)=\frac{p'(\rho)}{\rho}=\frac{c^2(\rho)}{\rho}
\end{equation*}
with $c(\rho)$ being the sound speed, and $p$ is the pressure that is a function of the density $\rho$.
For an ideal polytropic gas, the pressure $p$ and the sound speed $c$ are given by
$p(\rho)=\kappa \rho^{\gamma}$ and $c^2(\rho)=\kappa \gamma \rho^{\gamma-1}$
for constants $\gamma>1$ and $\kappa>0$.
Without loss of generality, we may choose $\kappa=1/\gamma$ to have
\begin{equation}
\label{1-c}
h(\rho)=\frac{\rho^{\gamma-1}-1}{\gamma-1},\quad c^2(\rho)=\rho^{\gamma-1}.
\end{equation}
This can be achieved by the following scaling:
$(t, {\bf x}, B)\longrightarrow (\alpha^2t, \alpha{\bf x},  \alpha^{-2}B)$ with $\alpha^2=\kappa\gamma$.
Taking the limit $\gamma\to 1+$, we can also consider the case of the isothermal flow ($\gamma=1$), for which
\begin{equation*}
i(\rho)=\ln \rho,\qquad c^2(\rho)\equiv 1.
\end{equation*}

In \S 2--\S 4, we first show how the shock problems can be formulated as free boundary problems
for the Euler equations for potential flow and discuss recently developed mathematical ideas, approaches,
and techniques for solving these free boundary problems.
Then, in \S 5, we present mathematical formulations of these shock problems for the
full Euler equations and discuss the role of the Euler equations for potential flow,
\eqref{1-a}--\eqref{1-b}, in these shock
problems in the realm of the full Euler equations.
Some further open problems in the direction
are also addressed.

\section{Shock Reflection-Diffraction and Free Boundary Problems}

\begin{figure}[h]
\centering
\includegraphics[height=0.9in,width=5.3in]{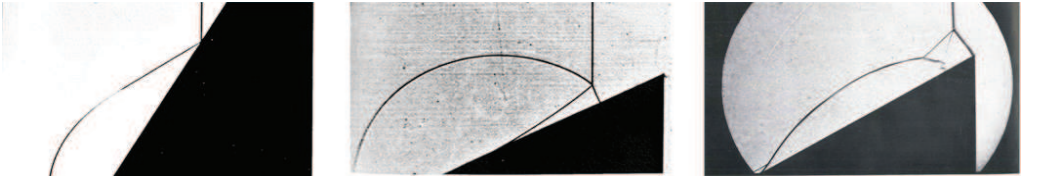}
\caption[]{Regular reflection (left); Simple Mach reflection (center); Irregular Mach reflection (right). From Van Dyke \cite{VD}, pp. 142--144}
 \label{fig:Experiment}
\end{figure}

We are first concerned with von Neumann's problem for
shock reflection-diffraction in \cite{Neumann1,Neumann2,Neumann}.
When a vertical planar shock perpendicular to the flow
direction $x_1$ and separating two uniform states (0) and (1),
with constant velocities
$(u_0, v_0)= (0, 0)$ and $(u_1,v_1)=(u_1, 0)$,
and constant densities $\rho_1>\rho_0$
(state (0) is ahead or to the right of the shock, and state
(1) is behind the shock), hits a symmetric wedge:
$$
W:= \{(x_1,x_2)\,:\, |x_2| < x_1 \tan \theta_{\rm w}, x_1 > 0\}
$$
head on at time $t = 0$,
a reflection-diffraction process takes place when $t > 0$.
Then a fundamental question is what type of wave patterns of
reflection-diffraction configurations may be formed around the wedge.
The complexity of reflection-diffraction configurations was first reported
by Ernst Mach \cite{Mach} in 1878, who first observed two patterns of
reflection-diffraction configurations:
Regular reflection (two-shock configuration; see Fig. \ref{fig:Experiment} (left))
and Mach reflection (three-shock/one-vortex-sheet configuration;
see Fig. \ref{fig:Experiment} (center)); also see \cite{BD,Chen-Feldman15,Courant-Friedrichs,VD}.
The issues remained dormant until the 1940s when John von Neumann \cite{Neumann1,Neumann2,Neumann},
as well as other mathematical/experimental
scientists ({\it cf.} \cite{BD,Chen-Feldman15,Courant-Friedrichs,Glimm-Majda,VD}
and the references cited therein)
began extensive research into all aspects of shock reflection-diffraction phenomena,
due to its importance in applications.
It has been found that the situations are much more complicated than
what Mach originally observed: The Mach reflection can be further
divided into more specific sub-patterns, and various other patterns of
shock reflection-diffraction configurations may occur such as the double Mach reflection,
von Neumann reflection, and Guderley reflection;
see \cite{BD,Chen-Feldman15,Courant-Friedrichs,Glimm-Majda,VD}
and the references cited therein.
Then the fundamental scientific issues include:
\begin{itemize}
\item[(i)] Structure of the shock reflection-diffraction configurations;

\item[(ii)] Transition criteria between the different patterns of shock
reflection-diffraction configurations;

\item[(iii)] Dependence of the patterns upon the physical parameters such as the
wedge-angle $\theta_{\rm w}$, the incident-shock-wave Mach number, and the
adiabatic exponent $\gamma\ge 1$.
\end{itemize}

\noindent
In particular, several transition criteria between the different
patterns of shock reflection-diffraction configurations have been proposed,
including the sonic conjecture and the detachment conjecture by von Neumann
\cite{Neumann1,Neumann2,Neumann}.

Careful asymptotic analysis has been made for various reflection-diffraction
configurations in Lighthill \cite{Lighthill-1,Lighthill-2},
Keller-Blank \cite{Keller-Blank}, Hunter-Keller \cite{Hunter-Keller}, Harabetian \cite{Harabetian},
Morawetz \cite{Morawetz}, and the references
cited therein; also see Glimm-Majda \cite{Glimm-Majda}.
Large or small scale numerical simulations have been also made;
{\it cf.} \cite{BD,Glimm-Majda,WC} and the references cited therein.
However, most of the fundamental issues for shock reflection-diffraction
phenomena have not been understood, especially the global structure and
transition between the different patterns of shock reflection-diffraction configurations.
This is partially because physical and numerical experiments are
hampered by various difficulties and have not yielded clear transition criteria
between the different patterns. In particular, numerical dissipation or physical
viscosity smear the shocks and cause boundary layers that interact with
the reflection-diffraction patterns and can cause spurious Mach steams; {\it cf.}
\cite{WC}. Furthermore, some different patterns occur when
the wedge angles are only fractions of a degree apart, a resolution even by
sophisticated modern experiments ({\it cf.} \cite{LD}) has not been able to reach.
For this reason, it is almost impossible to
distinguish experimentally between the sonic and detachment criteria,
as pointed out in \cite{BD}.
In this regard, the necessary approach to understand
fully the shock reflection-diffraction phenomena, especially the transition criteria,
is still via rigorous mathematical analysis.
To achieve this, it is essential to formulate the shock reflection-diffraction problem
as a free boundary problem and establish the global existence, regularity,
and structural stability of its solution.

\smallskip
Mathematically, the shock reflection-diffraction problem
is a multidimensional
lateral Riemann problem in domain
$\R^2\setminus \bar{W}$.

\begin{problemL}[Lateral Riemann Problem]\label{ibvp-c}
{\it
Piecewise constant initial data, consisting of state $(0)$
on $\{x_1>0\}\setminus \bar{W}$
and state $(1)$ on $\{x_1 < 0\}$ connected by a shock at $x_1=0$,
are prescribed at $t = 0$.
Seek a solution of the Euler system \eqref{1-a}--\eqref{1-b} for $t\ge 0$ subject to these initial
data and the boundary condition $\nabla\Phi\cdot\nnu=0$ on $\partial W$.
}
\end{problemL}

\medskip
Notice that {\bf Problem \ref{ibvp-c}}
is invariant under the scaling:
\begin{equation}\label{2.3a}
(t, {\bf x})\rightarrow (\alp t, \alp{\bf x}),\quad (\rho, \Phi)\rightarrow (\rho, \alpha \Phi)\qquad \text{for}\;\;\alp\neq 0.
\end{equation}
Thus, we seek self-similar solutions in the form of
\begin{equation}\label{1.2a}
\rho(t, {\bf x})=\rho(\xi,\eta),\quad \Phi(t, {\bf x})=t\phi(\xi,\eta)\qquad\text{for}\;\;(\xi,\eta)=\frac{{\bf x}}{t}.
\end{equation}
Then the pseudo-potential function $\vphi=\phi-\frac 12(\xi^2+\eta^2)$ satisfies
the following Euler equations for self-similar solutions:
\begin{align}
\label{1-r}
&{\rm div}(\rho D\vphi)+2\rho=0,\\
\label{1-p1}
&\frac{\rho^{\gam-1}-1}{\gam-1}+(\frac 12|D\vphi|^2+\vphi)=B,
\end{align}
where the divergence ${\rm div}$ and gradient $D$ are with respect to $(\xi,\eta)$.
From this, we obtain
the following second-order nonlinear PDE
for
$\vphi(\xi,\eta)$:
\begin{equation}
\label{2-1}
{\rm div}\big(\rho(|D\vphi|^2,\vphi)D\vphi\big)+2\rho(|D\vphi|^2,\vphi)=0
\end{equation}
with
\begin{equation}
\label{1-o}
\rho(|D\vphi|^2,\vphi)=
\bigl(B_0-(\gam-1)(\frac 12|D\vphi|^2+\vphi)\bigr)^{\frac{1}{\gam-1}},
\end{equation}
where $B_0:=(\gam-1)B+1$.
Then we have
\begin{equation}
\label{1-a1}
c^2(|D\vphi|^2,\vphi)=
B_0-(\gam-1)\big(\frac 12|D\vphi|^2+\vphi\big).
\end{equation}
Equation \eqref{2-1} is a {\it nonlinear PDE of mixed elliptic-hyperbolic type}.
It is elliptic if and only if
\begin{equation}
\label{1-f}
|D\vphi|<c(|D\vphi|^2,\vphi).
\end{equation}

If $\rho$ is a constant, then, by \eqref{2-1} and \eqref{1-o}, the corresponding
pseudo-potential $\vphi$ is in the form of
$$
\vphi(\xi,\eta)=-\frac 12(\xi^2+\eta^2)+u\xi+v\eta+k
$$
for constants $u,v$, and $k$.

\smallskip
Then {\bf Problem \ref{ibvp-c}} is reformulated as a boundary value problem
in unbounded domain
$$
\Lambda:=\R^2\setminus\{(\xi, \eta)\,:\,|\eta|\le \xi \tan\theta_{\rm w}, \xi>0\}
$$
in the self-similar coordinates $(\xi,\eta)$.

\begin{problemL}[Boundary Value Problem]\label{bvp-c}
{\it Seek a
solution $\varphi$ of equation \eqref{2-1} in the self-similar
domain $\Lambda$ with the slip boundary condition
$D\varphi\cdot\nnu|_{\partial\Lambda}=0$
on the wedge boundary $\partial\Lambda$
and the asymptotic boundary condition at infinity:
$$
\varphi\to\bar{\varphi}=
\begin{cases} \varphi_0 \qquad\mbox{for}\,\,\,
                         \xi>\xi_0, |\eta|>\xi \tan\theta_{\rm w},\\
              \varphi_1 \qquad \mbox{for}\,\,\,
                          \xi<\xi_0.
\end{cases}
\qquad \mbox{when $\xi^2+\eta^2\to \infty$,}
$$
where
\begin{equation*}
\pSi_0=-\frac{\xi^2+\eta^2}{2},
\qquad \pSi_1=-\frac{\xi^2+\eta^2}{2} +u_1(\xi-\xi_0),
\end{equation*}
and
$\xi_0=\rho_1\sqrt{\frac{2(c_1^2-c_0^2)}{(\gamma-1)(\rho_1^2-\rho_0^2)}}=\frac{\rho_1u_1}{\rho_1-\rho_0}$ is the location of the incident shock
in the $(\xi,\eta)$--coordinates.

By symmetry, we can restrict to the upper half-plane $\{\eta>0\}\cap\Lambda$,
with condition $\partial_\nnu \varphi=0$ on $\{\eta=0\}\cap\Lambda$.
}
\end{problemL}

A shock is a curve across which $D\vphi$ is discontinuous. If $\Om^+$ and $\Om^-(:=\Om\setminus \ol{\Om^+})$
are two nonempty open subsets of $\Om\subset \R^2$, and $S:=\der\Om^+\cap \Om$ is a $C^1$-curve
where $D\vphi$ has a jump, then $\vphi\in W^{1,1}_{\rm loc}\cap C^1(\Om^{\pm}\cup S)\cap C^2(\Om^{\pm})$
is a global weak solution of \eqref{2-1} in $\Om$ if and only if $\vphi$ is in $W^{1,\infty}_{loc}(\Om)$
and satisfies equation \eqref{2-1} and the Rankine-Hugoniot condition on $S$:
\begin{equation}
\label{1-h}
[\rho(|D\vphi|^2, \vphi)D\vphi\cdot\nnu_{\rm s}]_S=0,
\end{equation}
and the physical entropy condition:  {\it The density function $\rho(|D\vphi|^2, \vphi)$
increases across $S$ in the relative flow direction with respect to $S$},
where $[F]_S$ is defined by
$$
[F(\xi,\eta)]_S:=F(\xi,\eta)|_{\overline{\Om^-}}-F(\xi, \eta)|_{\overline{\Om^+}}\qquad\text{for}\;\;(\xi, \eta)\in S,
$$
and $\nnu_{\rm s}$ is a unit normal on $S$.

Note that the condition $\vphi\in W^{1,\infty}_{\rm loc}(\Om)$ requires
another Rankine-Hugoniot condition on $S$:
\begin{equation}
\label{1-i}
[\vphi]_S=0.
\end{equation}

\begin{figure}
\centering
\begin{minipage}{0.435\textwidth}
\centering
\psfrag{Om}{\Large$\Omega$}
\psfrag{P0}{$P_0$}
\psfrag{P1}{$P_1$}
\psfrag{P2}{$P_2$}
\psfrag{P3}{$P_3$}
\psfrag{P4}{$P_4$}
\includegraphics[height=1.3in,width=2.5in]{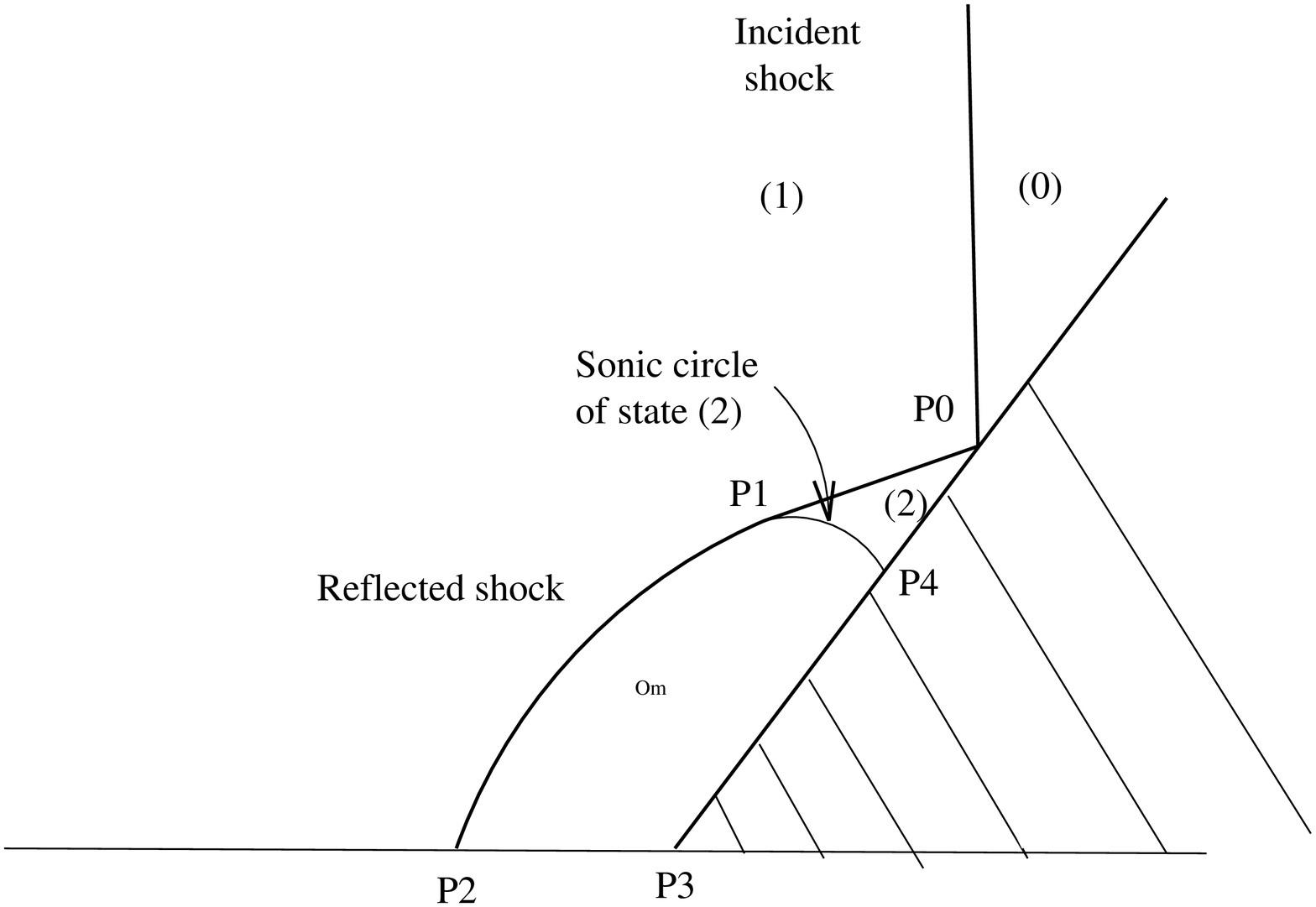}  %%% for latex (non-pdf)
\caption{Supersonic regular reflection}
\label{fig:RegularReflection}
\end{minipage}
\qquad\,\,\,
\begin{minipage}{0.435\textwidth}
\centering
\psfrag{Om}{\Large$\Omega$}
\psfrag{p0}{$P_0$}
\psfrag{p1}{$P_1$}
\psfrag{p2}{$P_2$}
\psfrag{p3}{$P_3$}
\psfrag{p4}{$P_4$}
\includegraphics[height=1.3in,width=2.4in]{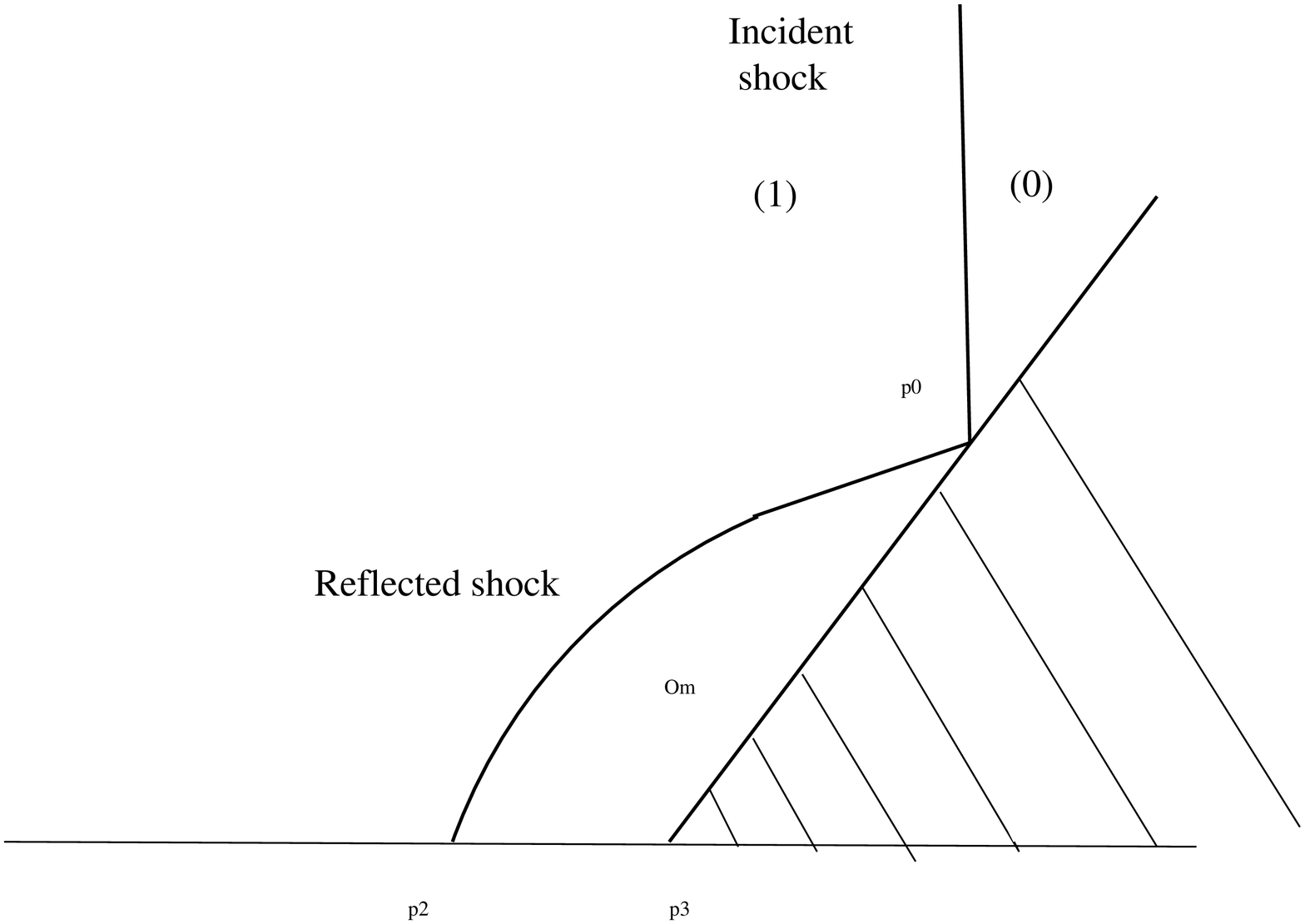}  %%% for latex (non-pdf)
\caption{Subsonic regular reflection}\label{fig:SubsoncRegularReflection}
\end{minipage}
\end{figure}

\smallskip
If a solution has one of the regular reflection-diffraction configurations as shown
in Figs. \ref{fig:RegularReflection}--\ref{fig:SubsoncRegularReflection},
and if $\varphi$ is smooth in the subregion between the wedge and reflected shock, then
it should satisfy  the boundary condition $D\varphi\cdot\nnu=0$ and the Rankine-Hugoniot
conditions \eqref{1-h}--\eqref{1-i} at $\PtIncW$
across the reflected shock separating it from state (1).
We define the uniform state (2) with pseudo-potential $\varphi_2(\xi,\eta)$ such
that
$$
\varphi_2(P_0)=\varphi_1(P_0), \quad D\varphi_2(\PtIncW)=D\varphi_1(\PtIncW),
$$
and the constant density $\rho_2$ of state (2) is equal to $\rho(|D\varphi|^2, \varphi)(\PtIncW)$
defined by (\ref{2-1}):
$$
\rho_2=\rho(|D\varphi|^2, \varphi)(\PtIncW).
$$
Then $D\varphi_2\cdot\nnu=0$ on the wedge boundary, and the Rankine-Hugoniot conditions
\eqref{1-h}--\eqref{1-i} hold on the flat shock $S_1=\{\varphi_1=\varphi_2\}$
between states (1) and (2), which passes through $\PtIncW$.

\smallskip
State (2) can be either subsonic or supersonic at $\PtIncW$.
This determines the subsonic or supersonic type of regular reflection-diffraction configurations.
The supersonic regular reflection-diffraction configuration
as shown in Fig. \ref{fig:RegularReflection}
consists of three uniform states (0), (1), (2), and a non-uniform state in domain $\Omega$,
where the equation is elliptic.
The reflected shock $\PtIncW\PtUpL\PtLwL$ has a straight part $\PtIncW\PtUpL$.
The elliptic domain $\Omega$ is separated from the hyperbolic region $\PtIncW\PtUpL\PtUpR$
of state (2) by a sonic arc $\PtUpL \PtUpR$.
The subsonic regular reflection-diffraction configuration
as shown in Fig. \ref{fig:SubsoncRegularReflection}
consists of two uniform states (0) and (1), and a non-uniform
state in domain $\Omega$, where the equation is elliptic, and
$\varphi_{|\Omega}(\PtIncW)=\varphi_2(\PtIncW)$ and
$D(\varphi_{|\Omega})(\PtIncW)=D\varphi_2(\PtIncW)$.

\smallskip
Thus, a necessary condition for the existence of regular reflection-diffraction
solution is the existence of the uniform state (2) determined by the conditions described above.
These conditions lead to an algebraic system for the constant velocity $(u_2, v_2)$ and
density $\rho_2$ of state (2),
which has solutions for some but not all of the wedge angles.
Specifically, for fixed densities $\rho_0<\rho_1$ of states (0) and (1),
there exist a sonic-angle $\theta_{\rm w}^{\rm s}$ and a detachment-angle
$\theta_{\rm w}^{\rm d}$
satisfying
$$
0<\theta_{\rm w}^{\rm d}<\theta_{\rm w}^{\rm s}<\frac{\pi}{2}
$$
such that state (2) exists for all $\theta_{\rm w}\in (\theta_{\rm w}^{\rm d}, \frac{\pi}{2})$
and does not exist for $\theta_{\rm w}\in (0, \theta_{\rm w}^{\rm d})$,
and the weak state (2) is supersonic at the reflection point $\PtIncW(\theta_{\rm w})$ for
$\theta_{\rm w}\in (\theta_{\rm w}^{\rm s}, \frac{\pi}{2})$,
sonic for $\theta_{\rm w}=\theta_{\rm w}^{\rm s}$,
and subsonic for $\theta_{\rm w}\in (\theta_{\rm w}^{\rm d}, \hat\theta_{\rm w}^{\rm s})$
for some $\hat\theta_{\rm w}^{\rm s}\in(\theta_{\rm w}^{\rm d}, \theta_{\rm w}^{\rm s}]$.

In fact, for each $\theta_{\rm w}\in(\theta_{\rm w}^{\rm d}, \frac{\pi}{2})$,
there exists also a {\em strong} state (2) with $\rho_2^{\rm strong}>\rho_2^{\rm weak}$.
There had been a long debate to
determine which one is physical for the local theory; see
\cite{BD,Courant-Friedrichs} and the references
cited therein.
It is expected that the strong reflection-diffraction configuration is non-physical;
indeed, it is shown as in Chen-Feldman \cite{Chen-Feldman10}
that the weak reflection-diffraction configuration tends to the unique normal
reflection, but the strong reflection-diffraction configuration does not,
when the wedge-angle $\theta_w$ tends
to $\frac{\pi}{2}$.
The strength of the corresponding reflected shock in the
weak reflection-diffraction configuration
is relatively weak
compared to the other shock given by the strong state (2), which is
called a {\it weak shock}.

\smallskip
If the weak state (2) is supersonic, the propagation speeds of the solution
are finite, and state (2) is
completely determined by the local information: State (1),
state (0), and the location of point $P_0$. That is, any information
from the region of reflection-diffraction,
especially the disturbance at corner $P_3$,
cannot travel towards the reflection point $P_0$.
However, if
it is subsonic, the information can reach $P_0$ and interact with
it, potentially altering a different reflection-diffraction configuration.
This argument motivated the following conjecture by
von Neumann in \cite{Neumann1,Neumann2}:

\medskip
{\bf The Sonic Conjecture}:
{\em There exists a supersonic reflection-diffraction
configuration when
$\theta_{\rm w}\in (\theta_{\rm w}^{\rm s}, \frac{\pi}{2})$
for $\theta_{\rm w}^{\rm s}>\theta_{\rm w}^{\rm d}$.
That is,
the supersonicity of the weak state {\rm (2)} implies the existence
of a supersonic regular reflection-diffraction
solution, as shown in Fig. {\rm \ref{fig:RegularReflection}}.}

\medskip
Another conjecture is that global regular reflection-diffraction
configuration is possible whenever the local regular reflection at the reflection
point is possible:

\medskip
{\bf The Detachment Conjecture}:
{\em There exists a regular reflection-diffraction configuration for
any wedge-angle $\theta_{\rm w}\in (\theta_{\rm w}^{\rm d}, \frac{\pi}{2})$.
That is, the existence of state {\rm (2)} implies the existence
of a regular reflection-diffraction
solution,
as shown in Figs. {\rm \ref{fig:RegularReflection}}--{\rm \ref{fig:SubsoncRegularReflection}}.
}

\medskip
It is clear that the supersonic/subsonic regular reflection-diffraction configurations are
not possible without a local two-shock configuration at the
reflection point on the wedge, so this is the weakest possible
criterion for the existence of supersonic/subsonic regular
shock reflection-diffraction configurations.

\begin{problemL}[Free Boundary Problem]\label{fbp-c}
{\it  For $\theta_{\rm w}\in (\theta_{\rm w}^d, \frac{\pi}{2})$,
find a free boundary (curved reflected shock) $P_1P_2$ on Fig. {\rm \ref{fig:RegularReflection}},
and $P_0P_2$ on Fig. {\rm \ref{fig:SubsoncRegularReflection}},
and a function $\varphi$ defined in region
$\Omega$ as shown in Figs. {\rm \ref{fig:RegularReflection}}--{\rm \ref{fig:SubsoncRegularReflection}},
such that $\varphi$ satisfies
\begin{itemize}
\item[\rm (i)]
Equation \eqref{2-1} in $\Om$;
\item[\rm (ii)]
$\vphi=\vphi_1$ and $\rho D\vphi\cdot\nnu_s=D\vphi_1\cdot\nnu_{\rm s}$ {on} the free boundary;
\item[\rm (iii)]
$\vphi=\vphi_2$ and $D\vphi=D\vphi_2$ {on} $P_1P_4$
in the supersonic case as shown in Fig. {\rm \ref{fig:RegularReflection}}
 and at $P_0$ in the subsonic case as shown in Fig. {\rm \ref{fig:SubsoncRegularReflection}};
\item[\rm (iv)]
$D\vphi\cdot\nnu=0$ {on} $\Wedge$,
\end{itemize}
where $\nnu_s$ and $\nnu$ are the interior unit normals to $\Omega$ on $\shock$ and $\Wedge$, respectively.
}
\end{problemL}

We observe that the key obstacle to the existence
of regular shock reflection-diffraction configurations as conjectured by von Neumann \cite{Neumann1,Neumann2}
is an additional possibility that,
for some wedge-angle $\theta_{\rm w}^{\rm a}\in (\theta_{\rm w}^{\rm d}, \frac{\pi}2)$, shock
$\PtIncW\PtLwL$ may attach to wedge-tip $\PtLwR$, as observed
by experimental results ({\it cf.} \cite[Fig. 238]{VD}).
To describe the conditions of such an attachment, we note that
$$
\rho_1>\rho_0, \qquad
u_1=(\rho_1-\rho_0)
\sqrt{\frac{2(\rho_1^{\gamma-1}-\rho_0^{\gamma-1})}{\rho_1^2-\rho_0^2}}.
$$
Then, for each $\rho_0$, there exists $\rho^{\rm c}>\rho_0$ such that
\begin{eqnarray*}
u_1\le c_1 \quad \mbox{if $\rho_1\in (\rho_0, \rho^{\rm c}]$}; \,\, \qquad
u_1>c_1 \quad \mbox{if $\rho_1\in (\rho^{\rm c}, \infty)$}.
\end{eqnarray*}

\smallskip
If $u_1\le c_1$, we can rule out the solution with a shock attached to the
wedge-tip.

\smallskip
If $u_1> c_1$, there would be a possibility that
the reflected shock could be attached to the wedge-tip
as experiments show ({\it e.g.} \cite[Fig. 238]{VD}).

\smallskip
\noindent
Thus, in \cite{Chen-Feldman10,Chen-Feldman15},
we have obtained the following results:

\smallskip
{\em
\begin{enumerate}
\item[\rm (i)] %\label{u1-le-c1}
If  $\rho_0$ and $\rho_1$ are such that $u_1\le c_1$, then the
supersonic/subsonic regular reflection-diffraction solution
exists for each wedge-angle $\theta_{\rm w}\in (\theta_{\rm w}^{\rm d}, \frac{\pi}{2})$;

\smallskip
\item[\rm (ii)]
If  $\rho_0$ and $\rho_1$ are such that $u_1> c_1$, then there exists
$\theta_{\rm w}^{\rm a}\in [\theta_{\rm w}^{\rm d}, \frac{\pi}2)$ such that
the regular reflection solution exists for
each  wedge-angle $\theta_{\rm w}\in (\theta_{\rm w}^{\rm a}, \frac{\pi}{2})$.
Moreover, if $\theta_{\rm w}^{\rm a}>\theta_{\rm w}^{\rm d}$,
then, for the wedge-angle $\theta_{\rm w}=\theta_{\rm w}^{\rm a}$,
there exists an {\it attached} solution, {\it i.e.},
a solution of {\rm \bf Problem \ref{fbp-c}} with  $\PtLwL=\PtLwR$.
\end{enumerate}
The type of regular reflection-diffraction configurations (supersonic as in Fig. {\rm \ref{fig:RegularReflection}},
or  subsonic as in Fig. {\rm \ref{fig:SubsoncRegularReflection}})
is determined by the type of state {\rm (2)} at $\PtIncW$.
For the supersonic and sonic reflection-diffraction case,
the reflected shock $\PtIncW\PtLwL$ is $C^{2,\alpha}$--smooth, and
the solution $\varphi$ is $C^{1,1}$ across the sonic arc for the supersonic
case,
which is optimal.
For the subsonic reflection-diffraction case
(Fig. {\rm \ref{fig:SubsoncRegularReflection}}),
the reflected shock $\PtIncW\PtLwL$
and the solution in $\Omega$ are both $C^{1,\alpha}$
near $P_0$ and $C^\infty$ away from $P_0$.
Furthermore, the regular reflection-diffraction solution tends to the unique normal
reflection, when wedge-angle $\theta_w$ tends to $\frac{\pi}{2}$.
}

\medskip
To solve this free boundary problem ({\bf Problem \ref{fbp-c}}),
we define a class of admissible solutions,
which are the solutions $\varphi$ with weak regular reflection-diffraction configurations,
such that, in the supersonic reflection case,
equation (\ref{2-1}) is strictly elliptic for $\varphi$
in $\overline\Omega\setminus \PtUpL\PtUpR$,
$\varphi_2\le\varphi\le \varphi_1$ holds in $\Omega$,
and the following monotonicity properties hold:
$$
\partial_\eta(\varphi_1-\varphi)\le 0, \quad D(\varphi_1-\varphi)\cdot {\bf e}\le 0 \qquad \mbox{in $\Omega\,\,\,$}
$$
for ${\bf e}=\frac{{P_0P_1}}{|P_0P_1|}$.
In the subsonic reflection case,
admissible solutions are defined similarly, with changes
corresponding to the structure of subsonic reflection-diffraction solution.

We derive uniform {\it a priori} estimates for admissible solutions with
any wedge-angle $\theta_{\rm w} \in [\theta_{\rm w}^{\rm d}+\varepsilon, \frac\pi 2]$
for each $\varepsilon>0$,
and then apply the degree theory to obtain
the existence for each $\theta_{\rm w} \in [\theta_{\rm w}^{\rm d}+\varepsilon, \frac\pi 2]$
in the class of admissible solutions,
starting from the unique normal reflection solution for $\theta_{\rm w}=\frac\pi 2$.
To derive the {\it a priori} bounds, we first obtain the estimates related to the geometry of the shock:
Show that the free boundary has a uniform positive distance from the sonic circle
of state (1) and from the wedge boundary away from $\PtLwL$ and $P_0$.
This allows to estimate the ellipticity of (\ref{2-1}) for $\varphi$ in $\Omega$
(depending on the distance to the sonic arc $P_1P_4$ for the supersonic reflection-diffraction configuration
and to $\PtIncW$ for the subsonic reflection-diffraction configuration).
Then we obtain the estimates near $P_1P_4$ (or $\PtIncW$
for the subsonic reflection) in scaled and weighted $C^{2,\alpha}$ for $\varphi$
and the free boundary, considering separately
four cases depending on $\frac{D\varphi_2}{c_2}$ at $\PtIncW$:

\begin{itemize}
\item[(i)] Supersonic:
$\frac{|D\varphi_2|}{c_2} \ge 1+\delta$;

\smallskip
\item[(ii)]
 Supersonic (almost sonic): $1<
\frac{|D\varphi_2|}{c_2} < 1+\delta$;

\smallskip
\item[(iii)]
Subsonic (almost sonic): $1-\delta\le \frac{ |D\varphi_2|}{c_2} \le 1$;

\smallskip
\item[(iv)]
Subsonic: $\frac{|D\varphi_2|}{c_2} \le 1-\delta$.
\end{itemize}

In cases (i)--(ii), equation (\ref{2-1}) is degenerate elliptic in $\Omega$
near  $\PtUpL\PtUpR$ on Fig. \ref{fig:RegularReflection}.
In case (iii), the equation is uniformly elliptic in $\overline\Omega$,
but the ellipticity constant is small near $\PtIncW$ on Fig. \ref{fig:SubsoncRegularReflection}.
Thus, in cases (i)--(iii), we use the local elliptic degeneracy,
which allows to find a comparison function in each case, to show
the appropriately fast decay of $\varphi-\varphi_2$ near $\PtUpL\PtUpR$ in cases (i)--(ii) and near $\PtIncW$ in case (iii);
furthermore, combining with appropriate local non-isotropic rescaling to obtain the uniform ellipticity,
we obtain the {\it a priori} estimates in the weighted and scaled $C^{2,\alpha}$--norms, which are different in each
of cases (i)--(iii),
but imply
the standard $C^{1,1}$--estimates in cases (i)--(ii), and the standard
$C^{2,\alpha}$--estimates in case (iii).
This is an extension of the methods developed in our earlier work \cite{Chen-Feldman10}.
In the uniformly elliptic case (iv), the solution is of subsonic reflection-diffraction
configuration as shown in Fig. \ref{fig:SubsoncRegularReflection},
and the estimates are more technically challenging than in cases (i)--(iii),
due to the lower {\it a priori} regularity
of the free boundary and since the uniform ellipticity does not allow
a comparison function that shows the decay of $\varphi-\varphi_2$ near $\PtIncW$.
Thus, we prove the $C^{\alpha}$--estimates of $D(\varphi-\varphi_2)$ near $P_0$.
With all of these, we provide a solution to the von Neumann's conjectures.

More details can be found in Chen-Feldman \cite{Chen-Feldman15};
also see \cite{BCF-09,Chen-Feldman10}.

\section{Shock Diffraction  (Lighthill's Problem) and Free Boundary Problems}

We are now concerned with shock diffraction by a two-dimensional wedge with convex corner
(Lighthill's problem).
When a plane shock in the $(t,\mathbf{x})$--coordinates,
$\mathbf{x}=(x_1,x_2)\in\mr^2$, with left state $(\r,u,v)=(\r_1,u_1,0)$
and right state $(\r_0,0,0)$ satisfying $u_1>0$ and $\r_0<\r_1$
from the left to right along the wedge with convex corner:
$$
W\deq\{(x_1,x_2)\, :\, x_2<0,\ x_1<x_2\tan\t_{\rm w}\},
$$
the incident shock interacts with the wedge
as it passes the corner, and then the shock diffraction occurs
({\it cf.} \cite{Lighthill-1,Lighthill-2}).
The mathematical study of the shock diffraction problem dates back to the 1950s by the work
of
Lighthill \cite{Lighthill-1,Lighthill-2} via asymptotic analysis;
also see \cite{Bargman,FWB,FTB}
via experimental analysis,
as well as
Courant-Friedrichs \cite{Courant-Friedrichs} and Whitham \cite{Whitham}.

\begin{figure}[!h]
  \centering
  \qquad \includegraphics[width=0.35\textwidth]{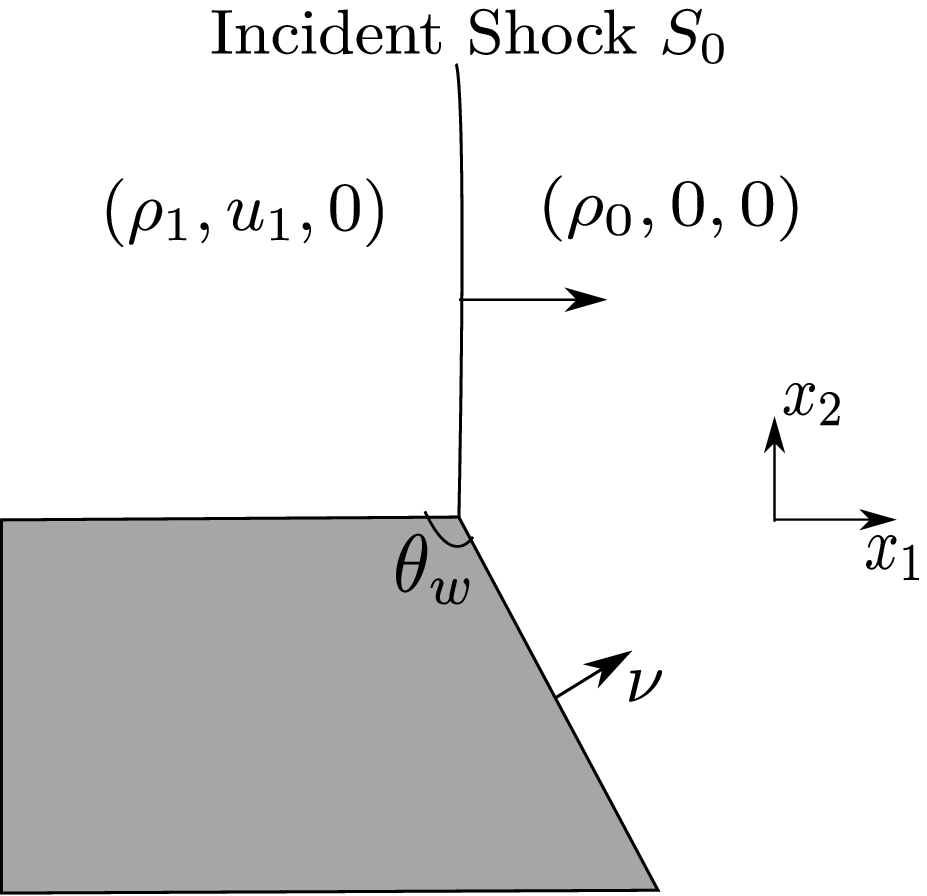}\qquad\qquad\qquad
  \includegraphics[width=0.35\textwidth]{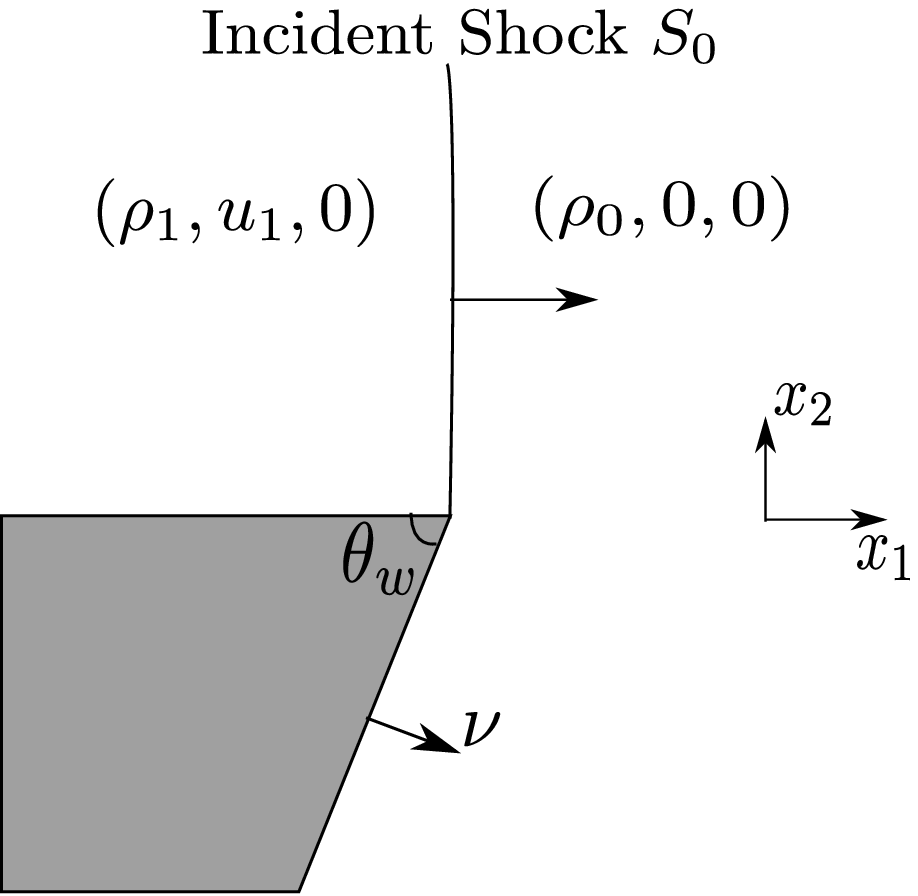}
  \caption{Lateral Riemann Problem}
\label{figure:initial}
\end{figure}

Similarly, this problem can be formulated as the
following lateral Riemann problem for potential flow:

\begin{problemL}[Lateral Riemann Problem; see Fig. \ref{figure:initial}] \label{ibvp-b}
Seek a solution of system \eqref{1-a}--\eqref{1-b}
with the initial condition at $t=0$:
\begin{equation}\label{con:1.1}
 (\r,\Phi)|_{t=0}=\begin{cases} (\r_1,u_1x_1)\qquad\,\, &\mbox{in}\,\, \{x_1<0, x_2>0\},\\[1mm]
(\r_0,0) \qquad\,\, &  \mbox{in}\,\,\{\theta_{\rm w}-\pi\leq\arctan\big(\frac{x_2}{x_1}\big)\leq\frac{\pi}{2}\},
\end{cases}
 \end{equation}
and the slip boundary condition along the wedge boundary $\partial W$:
 \begin{equation}\label{con:1.2}
 \nabla\Phi\cdot\nnu|_{\partial W}=0,
 \end{equation}
where $\nnu$ is the exterior unit normal to $\partial W$.
\end{problemL}

{\bf Problem \ref{ibvp-b}}
is also invariant under the self-similar scaling \eqref{2.3a}.
Thus, we seek self-similar solutions with form \eqref{1.2a}
in the self-similar domain outside the wedge:
$$
\L\deq\{\theta_w-\pi\leq\arctan\big(\frac{\eta}{\xi}\big)\leq\pi\}.
$$
Then the shock interacts with the pseudo-sonic circle of state $(1)$ to
become a transonic shock,
and {\bf Problem \ref{ibvp-b}}
can be formulated as the following
boundary value problem in the self-similar coordinates $(\xi,\eta)$.

\begin{figure}[!h]
  \centering
  \includegraphics[width=0.35\textwidth]{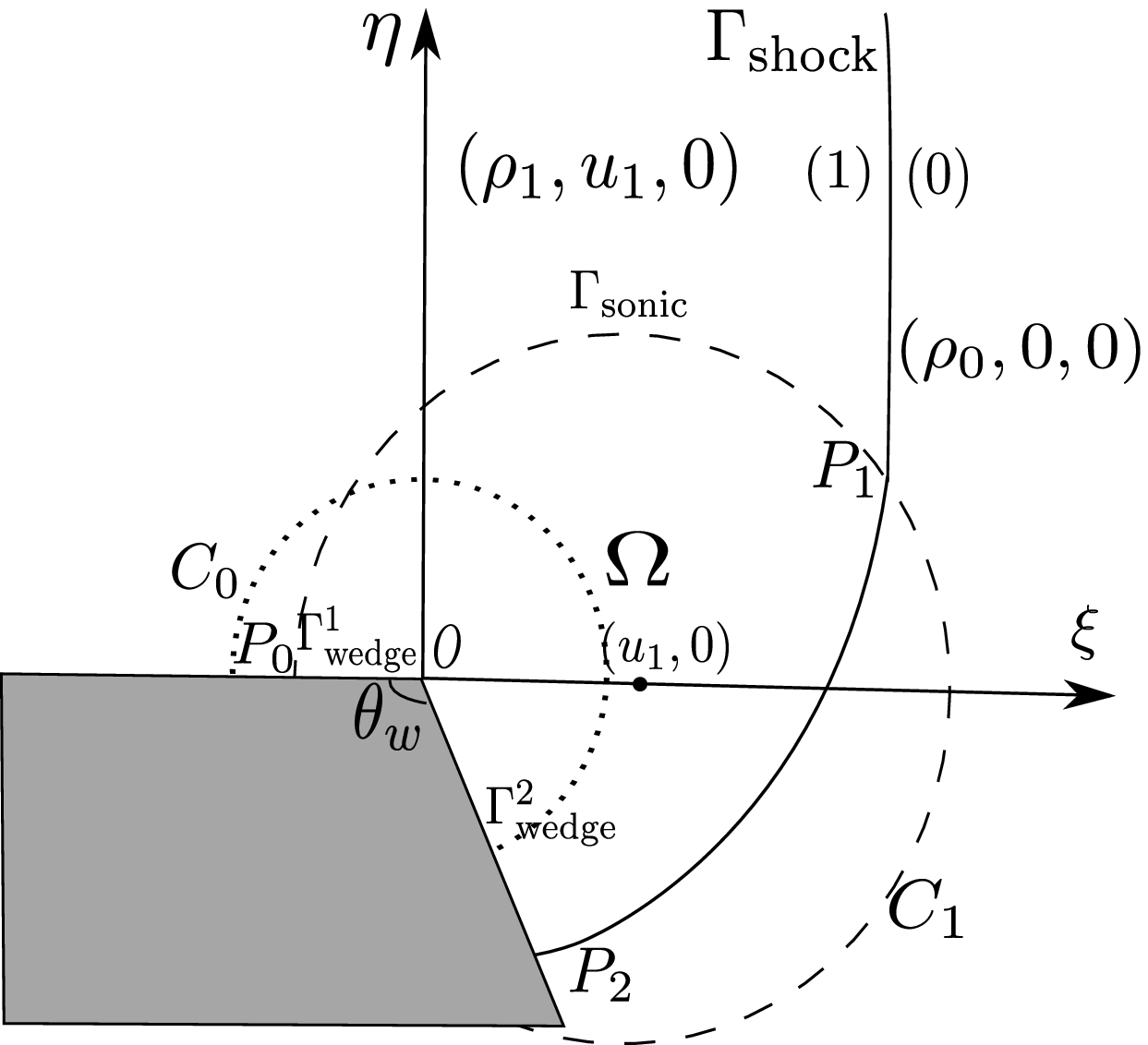}\qquad\qquad\qquad
  \includegraphics[width=0.35\textwidth]{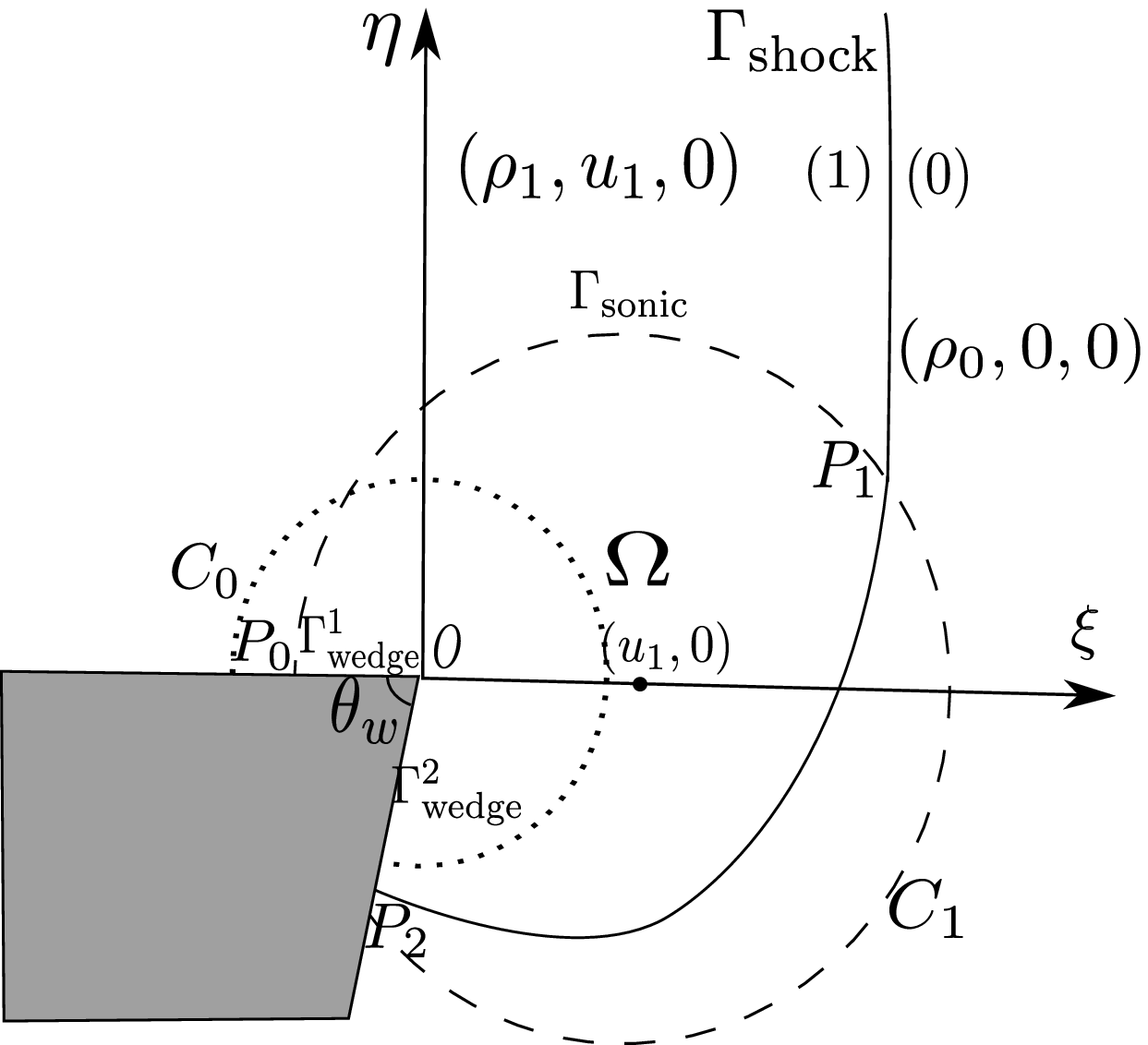}
  \caption{Boundary Value Problem; see Chen-Xiang \cite{Chen-Xiang-2}}\label{figure:interact}
\end{figure}

\begin{problemL}[Boundary Value Problem; see Fig. \ref{figure:interact}] \label{bvp-b}
Seek a solution $\pSi$ of equation \eqref{2-1}
in the self-similar
domain $\L$ with the
slip boundary condition on the wedge boundary $\partial\L$:
$$
D\pSi\cdot\nnu|_{\partial\L}=0
$$
and the asymptotic boundary condition at infinity:
$$
(\r,\pSi)\rightarrow (\bar{\r},\bar{\pSi})
=\begin{cases}
(\r_1, \pSi_1)\qquad&\text{in }\{\xi<\xi_0,\
\eta\geq0\},\\[1mm]
(\r_0, \pSi_0)\qquad&\text{in }\{\xi>\xi_0,\
\eta\geq0\}\cup\{\theta_{\rm w}-\pi\leq\arctan\big(\frac{\eta}{\xi}\big)\leq0\},
\end{cases}
$$
when $\xi^2+\eta^2\rightarrow \infty$ in the sense that
$
\lim_{R\rightarrow\infty}\|\pSi-\bar{\pSi}\|_{C^1(\L\backslash
B_{R}(0))}=0,
$
where $\pSi_0$, $\pSi_1$, and $\xi_0$ are the same as defined in {\rm {\bf Problem \ref{bvp-c}}}
in the $(\xi,\eta)$--coordinates.
\end{problemL}

Since $\pSi$ does not satisfy the slip boundary condition
for $\xi\geq0$, the solution must differ from state $(1)$ in
$\{\xi<\xi_1\}\cap\L$ near the wedge-corner, which forces the shock to be
diffracted by the wedge.
There is a critical angle $\theta_{c}$ so that, when $\theta_{\rm w}$ decreases to $\theta_{\rm c}$, two sonic circles $C_0$ and $C_1$
coincide at $P_2$ on $\Gamma_{\rm wedge}$.
Then {\bf Problem \ref{bvp-b}} can be formulated as
the following free boundary problem:

\begin{problemL}[Free Boundary Problem]\label{fbp-b}
For $\theta_{\rm w}\in (\theta_{\rm c}, \pi)$, find a free boundary (curved shock) $\shock$ and a function $\vphi$ defined in region $\Om$, enclosed
by $\shock, \Gamma_{\rm sonic}$, and the wedge boundary $\Gamma_{\rm wedge}:=\Gamma^1_{\rm wedge}\cup\Gamma^2_{\rm wedge}$,
such that $\vphi$ satisfies
\begin{itemize}
\item[\rm (i)]
Equation \eqref{2-1} in $\Om$;
\item[\rm (ii)]
$\vphi=\vphi_0$, $\rho D\vphi\cdot\nnu_{\rm s}=\rho_0D\vphi_0\cdot\nnu_{\rm s}$ {on} $\shock$;
\item[\rm (iii)]
$\vphi=\vphi_1$, $D\vphi=D\vphi_1$ {on} $\Gamma_{\rm sonic}$;
\item[\rm (iv)]
$D\vphi\cdot\nnu=0$ {on} $\Wedge$,
\end{itemize}
where $\nnu_s$ and $\nnu$ are the interior unit normals to $\Omega$ on $\shock$ and $\Wedge$, respectively.
\end{problemL}

In domain $\O$,
the solution is
expected to be pseudo-subsonic and smooth, to satisfy the slip boundary
condition along the wedge, and to be $C^{1,1}$--continuous across
the pseudo-sonic circle to become pseudo-supersonic.
Then the solution of {\bf Problem \ref{fbp-b}} can be shown to be
the solution of {\bf Problem \ref{ibvp-b}}.

The free boundary problem has been solved in \cite{Chen-Xiang-1,Chen-Xiang-2}.
A crucial challenge of this problem is that
the expected elliptic domain of the solution is concave so that its boundary
does not satisfy the exterior ball condition, since the angle $2\pi-\theta_w$
exterior to the wedge at the origin is larger than $\pi$ for the given wedge-angle
$\theta_w\in (0, \pi)$,
besides other mathematical difficulties
including free boundary problems without uniform oblique derivative conditions.
There is no general theory of elliptic PDEs
on such concave domains, whose coefficients involve the gradient of
the solutions.
In general, the expected regularity in this domain,
even for Laplace's equation,
is only $C^{\a}$ with $\a<1$;
however, the coefficients in \eqref{2-1} depend on the gradient of $\pSi$
so that the ellipticity of
this equation depends also on the boundedness of the derivatives,
which is one of the essential
difficulties of this problem.
To overcome the difficulty, the physical boundary conditions must
be exploited to force a finer regularity of solutions at the corner to
let equation \eqref{2-1} make sense.
More precisely, the
strategy here is that, instead of analyzing
equation \eqref{2-1} directly, we study another system of equations
for the physical quantities $(\r, u, v)$
for the existence of the velocity potential.

A tempting try would be to differentiate first
equation \eqref{1-r}
to obtain an equation for
$v$, then use the irrotationality to solve $u$ (once $v$ has solved),
and finally use \eqref{1-p1}
to solve the density $\r$.
In order to show the equivalence between these equations
and the original potential flow equation \eqref{2-1},
an additional one-point boundary condition is required for $v$.
However, it is unclear for the boundary condition to be deduced for $v$
for the problem.
Moreover, along the boundaries $\Gsh$ and $\Gwt$
which meet at the corner,
the derivative boundary conditions of the deduced
second-order elliptic equation to $v$ are the second kind boundary
conditions, {\it i.e.} without the viscosity, compared to
\cite{Chen-Feldman10}.
This implies that the results from
\cite{l-advmath1985161Perronobliqueboundary, l-jmathanalappl1986422Perronmixedboundary}
could not be directly used.
To overcome this,
the following directional velocity $(w,z)$ is introduced
whose relation with $(u,v)$ is
\begin{align*}
(w, z):=(u\sin\theta_{\rm w}-v\cos\theta_{\rm w}, u\cos\theta_{\rm w}+v\sin\theta_{\rm w}),
\end{align*}
such that the one-point boundary condition for $w$
is not required for solving $w$, and then treat $z$ as $u$.
For $(w,z)$, the $C^{\a}$--regularity is enough.

On the other hand, for these equations, some new technical difficulties arise,
for which
new mathematical ideas and techniques have to be developed.
First, it is a coupled system so that
the coefficients of the nonlinear degenerate elliptic equation for $w$ depend on $z$,
which makes the uniform estimates for $w$ near the sonic circle more challenging.
Second, the obliqueness condition on the free boundary deduced from
the Rankine-Hugoniot conditions depends on the smallness of $z$.
To overcome this,
a degenerate elliptic cut-off function near the pseudo--sonic circle
is introduced,
which is more precise in comparison with
\cite{Chen-Feldman10}.
The reason why the more precise degenerate cut-off function
requires to be introduced
is that the uniform estimates of $w$ are required
to obtain a convergent
sequence near $P_1$, which is crucial for
the equivalence between the deduced system and the potential flow
equation \eqref{2-1} with degenerate elliptic cut-off.
Third, since the new feature that $\sin\theta$ may be $0$
along the pseudo--sonic circle and the fact that there
is no $C^2$--regularity at $P_1$ where the shock and pseudo--sonic
circle meet from the optimal regularity argument
by Bae-Chen-Feldman \cite{BCF-09},
more effort is needed to remove the degenerate
elliptic cut-off case by case carefully,
near and away from $P_1$ respectively.
The final main difficulty
is to show the equivalence of the original potential flow equation \eqref{2-1}
and the deduced second-order equation for $w$ with irrotationality and
Bernoulli's law, which requires gradient estimates for $w$ near
the pseudo-sonic circle, but the estimates by scaling only provide
a bound divided by the distance to this circle.
This is overcome thanks to the estimates involving $\e$.

When the wedge-angle becomes smaller,
several other difficulties arise.
Due to the concave corner at the origin,
more technical arguments are required to obtain the existence
of solutions to the modified problem.
Unlike in \cite{Chen-Feldman15},
since it requires to take derivatives along the shock to obtain
a boundary condition for $w$,
a new way to modify the Rankine-Hugoniot conditions
is designed delicately,
based on the nonlinear structure of the shock.
From this modified condition, the Dirichlet condition is assigned on the shock
where the modified uniform oblique condition fails.
Thus,
the uniform boundedness of solutions need to be controlled more carefully.
Finally,
the existence of shock diffraction configuration
up to the critical angle $\theta_{\rm c}$  is established.

\section{Prandtl-Meyer Reflection Configurations and Free Boundary Problems}

We now consider with Prandtl-Meyer's problem for unsteady global solutions for supersonic flow past a solid ramp,
which can be also regarded as portraying the symmetric gas flow impinging onto a solid wedge (by symmetry).
When a steady supersonic flow past a solid ramp whose angle is less than
the critical angle (called the detachment angle) $\theta_{\rm w}^d$,
Prandtl \cite{Prandtl} employed the shock polar analysis to show that there are two possible
configurations:
The weak shock reflection with supersonic or subsonic downstream flow
and the strong shock reflection
with subsonic downstream flow, which both satisfy the physical entropy
conditions, provided that we do not give additional conditions at downstream;
also see \cite{Busemann,Courant-Friedrichs,Meyer}.

The fundamental question of whether one or both of the strong and the weak shocks are physically
admissible has been vigorously debated over the past seventy years, but has not yet been settled
in a definite manner ({\it cf.} \cite{Courant-Friedrichs,Dafermos10,Serre}).
On the basis of experimental and numerical evidence, there are strong indications that it is
the weak reflection that is physically admissible.
One plausible approach is to single out the strong shock reflection by the consideration
of stability: The stable ones are physical.
It has been shown in the steady regime that the weak reflection
is not only structurally stable ({\it cf.} \cite{CZZ}),
but also $L^1$-stable with respect to steady small perturbation of both the ramp slope
and the incoming steady upstream flow ({\it cf.} \cite{CLi}),
while the strong reflection is also structurally stable for a large spectrum of
physical parameters ({\it cf.} \cite{Chen-Chen-Feldman,C-Fang}).

We are interested in the rigorous unsteady analysis of the steady supersonic weak shock solution
as the long-time behavior of an unsteady flow
and establishing the stability of the steady supersonic weak shock solution as
the long-time asymptotics of an unsteady flow with the Prandtl-Meyer configuration
for all the admissible physical parameters for potential flow.
Our goal is to find a solution $(\rho, \Phi)$ to system \eqref{1-a}--\eqref{1-b}
when a uniform flow in $\R^2_+:=\{x_1\in \R, x_2>0\}$ with
$(\rho, \nabla_{\bf x}\Phi)=(\irho, \iu,0), \x=(x_1, x_2)$, is heading to a solid ramp
at $t=0$:
$$
W:=\{(x_1,x_2)\,:\, 0<x_2<x_1\tan\theta_{\rm w}, x_1>0\}.
$$

\begin{problemL}
[Lateral Riemann Problem]
\label{problem-1}
Seek a solution of system \eqref{1-a}--\eqref{1-b} with
$B=\frac{\iu^2}{2}+\frac{\irho^{\gam-1}-1}{\gam-1}$ and the initial condition at $t=0$:
\begin{equation}
\label{1-d}
(\rho,\Phi)|_{t=0}=(\irho, \iu x_1)
\qquad\text{for}\;\;(x_1,x_2)\in \R^2_+\setminus W,
\end{equation}
and with the slip boundary condition along the wedge boundary $\der W$:
\begin{equation}
\label{1-e}
\nabla_{\bf x}\Phi\cdot \nnu|_{\der W\cap\{x_2>0\}}=0,
\end{equation}
where $\nnu$ is the exterior unit normal to $\der W$.
\end{problemL}

Again, {\bf Problem \ref{problem-1}} is
is invariant under the self-similar scaling \eqref{2.3a}.
Thus, we seek self-similar solutions in form \eqref{1.2a}
so that the pseudo-potential function $\vphi=\phi-\frac 12(\xi^2+\eta^2)$ satisfies
the nonlinear PDE \eqref{2-1} of mixed type.

As the incoming flow has the constant velocity $(\iu,0)$, the corresponding pseudo-potential
$\ivphi$ has the expression of
\begin{equation}
\label{1-m}
\ivphi=-\frac 12(\xi^2+\eta^2)+\iu\xi.
%+\ik
\end{equation}

Then {\bf Problem \ref{problem-1}} can be reformulated as the following boundary value problem
in the domain
$$
\Lambda:=\R^2_+\setminus\{(\xi,\eta)\,: \,\eta\le \xi\tan\theta_{\rm w},\, \xi\ge 0\}.
$$
in the
self-similar coordinates $(\xi,\eta)$, which corresponds to $\{(t, {\bf x})\, :\, {\bf x}\in \R^2_+\setminus W,\, t>0\}$
in the $(t, {\bf x})$--coordinates:

\medskip
\begin{problemL}[Boundary Value Problem]
\label{problem-2}
Seek a solution $\vphi$ of equation \eqref{2-1} in the self-similar domain $\Lambda$
with the slip boundary condition:
\begin{equation}
\label{1-k}
D\vphi\cdot\nnu|_{\partial \Lambda}=0,
\end{equation}
and the asymptotic boundary condition at infinity:
\begin{equation}\label{1-k-b}
\vphi-\vphi_\infty\longrightarrow 0
\end{equation}
along each ray $R_\theta:=\{ \xi=\eta \cot \theta, \eta > 0 \}$
with $\theta\in (\theta_{\rm w}, \pi)$ as $\eta\to \infty$
in the sense that
\begin{equation}\label{1-k-c}
\lim_{r\to \infty} \|\varphi  - \varphi_\infty \|_{C(R_{\theta}\setminus B_r(0))} = 0.
\end{equation}
\end{problemL}

\begin{figure}[htp]
\centering
\begin{psfrags}
\psfrag{ls}[cc][][0.8][0]{$\leftshockop$}
\psfrag{sn}[cc][][0.8][0]{$\rightshockop$}
\psfrag{ol}[cc][][0.8][0]{$\oOmop$}
\psfrag{on}[cc][][0.8][0]{$\nOmop$}
\psfrag{lsn}[cc][][0.8][0]{$\phantom{aa}\leftsonicop$}
\psfrag{rsn}[cc][][0.8][0]{$\rightsonicop\phantom{aaaaaaa}$}
\psfrag{tw}[cc][][0.8][0]{$\theta_w$}
\psfrag{o}[cc][][0.8][0]{$O$}
\psfrag{lb}[cc][][0.8][0]{$\leftbottom$}
\psfrag{rb}[cc][][0.8][0]{$\rightbottom$}
\psfrag{lt}[cc][][0.8][0]{$\lefttop$}
\psfrag{rt}[cc][][0.8][0]{$\righttop$}
\psfrag{s}[cc][][0.8][0]{$\shock$}
\psfrag{om}[cc][][0.8][0]{$\Om$}
\includegraphics[scale=.5]{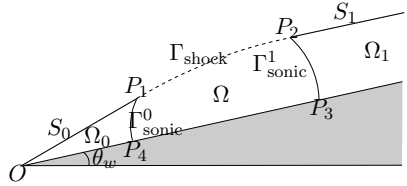}
\caption{Admissible solutions for $\theta_{\rm w}\in (0, \theta_{\rm w}^s)$ in the self-similar coordinates $(\xi, \eta)$;
see Bae-Chen-Feldman \cite{BCF-14}}\label{fig:global-structure-1}
\end{psfrags}
\end{figure}

\begin{figure}[htp]
\centering
\begin{psfrags}
\psfrag{oi}[cc][][0.8][0]{$\phantom{aaaaaaa}\Om_{\infty}: (\iu, 0), \rho_\infty$}
\psfrag{om}[cc][][0.8][0]{$\Om$}
\psfrag{on}[cc][][0.8][0]{$\nOmop$}
\psfrag{sn}[cc][][0.8][0]{$\rightshockop$}
\psfrag{ls}[cc][][0.8][0]{$\leftshockop$}
\psfrag{s}[cc][][0.8][0]{\phantom{aaaaa}$\shock$}
\psfrag{tw}[cc][][0.8][0]{$\theta_w$}
\includegraphics[scale=0.7]{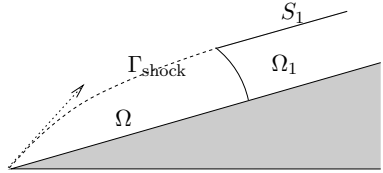}
\caption{Admissible solutions for $\theta_{\rm w}\in [\theta_{\rm w}^s,\theta_{\rm w}^d)$ in the self-similar coordinates $(\xi, \eta)$;
see Bae-Chen-Feldman \cite{BCF-14}}\label{fig:global-structure-2}
%\caption{Weak shock solutions in $(\xi,\eta)$-coordinates($\theta_{sonic}\le \theta_w<\theta_{detach}$)}\label{fig:global-structure-2}
%\caption{}\label{}
\end{psfrags}
\end{figure}

In particular, we seek a weak solution of {\bf Problem \ref{problem-2}}
with two types of Prandtl-Meyer reflection configurations whose occurrence
is determined by wedge-angle $\theta_{\rm w}$ for the two different cases:
One contains a straight weak oblique shock attached to wedge-tip $O$ and
the oblique shock is connected to a normal shock through a curved shock
when $\theta_{\rm w}<\theta_{\rm w}^s$,
as shown in Fig. \ref{fig:global-structure-1}; the other contains a curved shock attached
to the wedge-tip
and connected to a normal shock when $\theta_{\rm w}^s\le \theta_{\rm w}<\theta_{\rm w}^d$,
as shown in Fig. \ref{fig:global-structure-2}, in which
the curved shock $\Gamma_{\rm shock}$ is tangential to
a straight weak oblique shock $S_0$ at the wedge-tip.

To seek a global entropy solution of {\bf Problem \ref{problem-2}} with the structure
of Fig. \ref{fig:global-structure-1} or Fig. \ref{fig:global-structure-2},
one needs to compute the pseudo-potential $\vphi_0$ below $S_0$.

\begin{figure}[htp]
\centering
\begin{psfrags}
\psfrag{tc}[cc][][0.8][0]{$\phantom{aaaaaaaaaaa}v=u\tan\theta_{\rm w}^s$}
\psfrag{tw}[cc][][0.8][0]{$\phantom{aaaaaaaaaaa}v=u \tan\theta_w$}
\psfrag{td}[cc][][0.8][0]{$\phantom{aaaaaaaaaaa}v=u\tan\theta_{\rm w}^d$}
\psfrag{u}[cc][][0.8][0]{$u$}
\psfrag{v}[cc][][0.8][0]{$v$}
\psfrag{u0}[cc][][0.8][0]{$\iu$}
\psfrag{zeta}[cc][][0.8][0]{}
\psfrag{ud}[cc][][0.8][0]{$u_{\rm d}$}
\includegraphics[scale=.8]{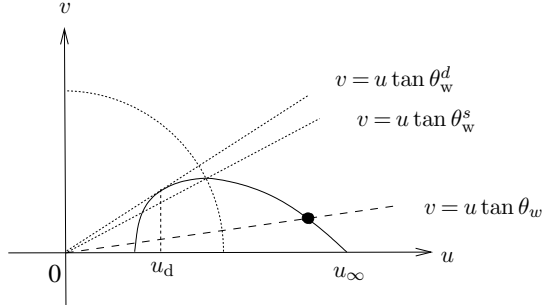}
\caption{Shock polars in the $(u,v)$--plane}\label{fig:polar}
\end{psfrags}
\end{figure}

Given $M_\infty>1$, we obtain $(u_0,v_0)$ and $\rho_0$ by using the shock polar curve
in Fig. \ref{fig:polar} for steady potential flow.
In Fig. \ref{fig:polar},
$\theta_{\rm w}^s$ is the wedge-angle such that
line $v=u \tan\theta_{\rm w}^s$ intersects
with the shock polar curve at a point on the circle of radius $c_\infty$,
and $\theta_{\rm w}^d$ is the wedge-angle so that
line $v=u \tan\theta_{\rm w}^d$
is tangential to the shock polar curve and there is no intersection between
line $v=u \tan \theta_{\rm w}$ and
the shock polar when $\theta_{\rm w}>\theta_{\rm w}^d$.

For any wedge-angle $\theta_{\rm w}\in (0,\theta_{\rm w}^s)$,
line $v=u\tan\theta_{\rm w}$ and the shock polar
curve intersect at a point $(u_0, v_0)$ with $|(u_0,v_0)|
>c_\infty$ and $u_0<\iu$;
while, for any wedge-angle $\theta_{\rm w}\in [\theta_{\rm w}^s, \theta_{\rm w}^d)$,
they intersect at a point $(u_0, v_0)$ with $u_0>u_{\rm w}^d$
and $|(u_0,v_0)|<c_\infty$.
The intersection state $(u_0, v_0)$ is the velocity for steady potential flow
behind an oblique shock $S_0$ attached to the wedge-tip with angle $\theta_{\rm w}$.
The strength of shock $S_0$ is relatively weak compared to the other shock
given by the other intersection point on the shock polar curve, which is a
a \emph{weak shock},
and the corresponding state $(u_0, v_0)$ is a \emph{weak state}.

We also note that
states
$(u_0, v_0)$ smoothly depend on $\iu$ and $\theta_{\rm w}$, and
such states are supersonic when $\theta_{\rm w}\in (0,\theta_{\rm w}^s)$ and
subsonic when $\theta_{\rm w}\in [\theta_{\rm w}^s, \theta_{\rm w}^d)$.

Once $(u_0, v_0)$ is determined,
by \eqref{1-i} and \eqref{1-m},
the pseudo-potentials $\vphi_0$ and $\vphi_1$ below the weak oblique shock $S_0$ and the normal shock $S_1$
are respectively in the form of
\begin{equation}
\label{1-n}
\vphi_0=-\frac 12(\xi^2+\eta^2)+u_0\xi+ v_0\eta,\quad
\vphi_1=-\frac 12(\xi^2+\eta^2)+ u_1 \xi+v_1\eta+ k_1
\end{equation}
for constants $u_0, v_0, u_1, v_1$, and $k_1$. Then it follows
from \eqref{1-o} and \eqref{1-n} that the corresponding densities $\rho_0$ and $\rho_1$
below $S_0$ and $S_1$ are constants, respectively.
In particular, we have
\begin{equation}
\label{2-n1}
\rho_0^{\gam-1}=\rho_\infty^{\gamma-1}+\frac{\gam-1}{2}\big(\iu^2-u_0^2-v_0^2\big).
\end{equation}

Then {\bf Problem \ref{problem-2}} can be formulated into the following
free boundary problem.

\smallskip
\begin{problemL}[Free Boundary Problem]\label{fbp-a}
For $\theta_{\rm w}\in (0, \theta_{\rm w}^d)$,
find a free boundary (curved shock) $\shock$ and a function $\vphi$ defined in domain $\Om$,
as shown in Figs. \ref{fig:global-structure-1}--\ref{fig:global-structure-2},
%enclosed by $\shock, \leftsonic, \rightsonic$, and $\Wedge$,
such that $\vphi$ satisfies
\begin{itemize}
\item[\rm (i)]
Equation \eqref{2-1} in $\Om$;
\item[\rm (ii)]
$\vphi=\ivphi$ and $\rho D\vphi\cdot\nnu_{\rm s}=D\ivphi\cdot\nnu_{\rm s}$ {on} $\shock$;
\item[\rm (iii)]
$\vphi=\hat{\vphi}$ and $D\vphi=D\hat{\vphi}$ {on} $\Gamma_{\rm sonic}^0\cup\Gamma_{\rm sonic}^1\,\,$
when $\theta_{\rm w}\in (0, \theta_{\rm w}^s)$
and on $\Gamma_{\rm sonic}^1\cup O\,\,$ when $\theta_{\rm w}\in [\theta_{\rm w}^s, \theta_{\rm w}^d)$
for $\hat{\vphi}:=\max(\vphi_0, \vphi_1)$;
\item[\rm (iv)]
$D\vphi\cdot \nnu=0$ {on} $\Wedge$,
\end{itemize}
where $\nnu_{\rm s}$ and $\nnu$ are the interior unit normals to $\Omega$ on $\shock$ and $\Wedge$, respectively.
\end{problemL}

Let $\vphi$ be a solution of {\bf Problem \ref{fbp-a}} with shock $\shock$.
Moreover, assume that $\vphi\in C^1(\overline\Om)$,
and $\shock$ is a $C^1$--curve up to its endpoints.
To obtain a solution of {\bf Problem \ref{problem-2}} from $\vphi$,
we have two cases:

For $\theta_{\rm w}\in (0, \theta_{\rm w}^s)$,
we divide half-plane $\{\eta\ge 0\}$ into five separate regions.
Let $\Om_S$ be the unbounded domain below curve $\ol{S_0\cup\shock\cup S_1}$ and above $\Wedge$
%in $\{\eta\ge 0\}$
(see Fig. \ref{fig:global-structure-1}).
In $\Om_S$, let $\Om_0$ be the bounded open domain enclosed by $S_0, \Gamma^0_{\rm sonic}$,
and $\{\eta=0\}$.
We set $\Om_{1}:=\Om_S\setminus \ol{(\Om_0\cup\Om)}$.
Define a function $\vphi_*$ in $\{\eta\ge 0\}$ by
\begin{equation}\label{extsol}
\vphi_*=
\begin{cases}
\ivphi& \qquad \text{in}\,  \Lambda\cap\{\eta\ge 0\}\setminus \Om_S,\\
\vphi_0& \qquad \text{in}\,\Om_0,\\
\vphi& \qquad \text{in}\, \Gamma^0_{\rm sonic}\cup\Om\cup\Gamma^1_{\rm sonic},\\
\vphi_1&\qquad \text{in}\,\Om_1.
\end{cases}
\end{equation}
By \eqref{1-i} and (iii) of {\bf Problem \ref{fbp-a}}, $\vphi_*$ is continuous in $\{\eta\ge 0\}\setminus\Om_S$
and is $C^1$ in $\overline{\Om_S}$.
In particular, $\vphi_*$ is $C^1$ across $\Gamma^0_{\rm sonic}\cup\Gamma^1_{\rm sonic}$.
Moreover, using (i)--(iii) of {\bf Problem \ref{fbp-a}}, we obtain
that  $\vphi_*$ is a global entropy solution of equation \eqref{2-1}
in $\Lambda\cap \{\eta>0\}$, which is the
Prandtl-Meyer's supersonic reflection configuration.

For $\theta_w\in [\theta_{\rm w}^s, \theta_{\rm w}^d)$,
region $\Omega_0\cup \Gamma^0_{\rm sonic}$ in $\varphi_*$ reduces to one point $O$,
and the corresponding $\varphi_*$ is a global entropy solution of equation \eqref{2-1}
in $\Lambda\cap \{\eta>0\}$,
which is the
Prandtl-Meyer's subsonic reflection configuration.

\medskip
The free boundary problem ({\bf Problem  \ref{fbp-a}}) has been solved in Bae-Chen-Feldman \cite{BCF-13,BCF-14}.

\smallskip
To solve this free boundary problem,
we follow the approach introduced in
Chen-Feldman \cite{Chen-Feldman15}.
We first define a class of admissible solutions,
which are the solutions $\varphi$ with Prandtl-Meyer reflection configuration,
such that,
when $\theta_{\rm w}\in (0, \theta_{\rm w}^s)$,
equation (\ref{2-1}) is strictly elliptic for $\varphi$
in $\overline\Omega\setminus (\Gamma_{\rm sonic}^{0}\cup\Gamma_{\rm sonic}^{1})$,
$\max\{\varphi_{0}, \varphi_{1}\}\le\varphi\le \varphi_\infty$ holds in $\Omega$,
and the following monotonicity properties hold:
$$
D(\varphi_\infty-\varphi)\cdot \mathbf{e}_{S_1}\ge 0, \quad D(\varphi_\infty-\varphi)\cdot \mathbf{e}_{S_0}\le 0 \qquad \mbox{in $\Omega$},
$$
where $\mathbf{e}_{S_0}$ and $\mathbf{e}_{S_{1}}$ are the unit tangential
directions to lines $S_0$ and $S_1$, respectively,
pointing to the positive $\xi$-direction.
For the case $\theta_{\rm w}\in [\theta_{\rm w}^s, \theta_{\rm w}^d)$,
admissible solutions are defined similarly, with
corresponding changes to the structure of subsonic reflection solutions.

We derive uniform {\it a priori} estimates for admissible solutions
for any wedge-angle $\theta_{\rm w} \in [0, \theta_{\rm w}^d-\varepsilon]$
for each $\varepsilon>0$,
and then employ the Leray-Schauder degree argument to obtain
the existence for each $\theta_{\rm w} \in [0, \theta_{\rm w}^d-\varepsilon]$
in the class of admissible solutions,
starting from the unique normal solution for $\theta_{\rm w}=0$.

More details can be found in \cite{BCF-13,BCF-14}; also see \S 2 above and
Chen-Feldman \cite{Chen-Feldman15}.

\medskip
In Chen-Feldman-Xiang \cite{Chen-Feldman-Xiang},
we have also established the strict convexity of the curved (transonic)
part of the free boundary in
the shock reflection-diffraction problem in \S 2 (also see Chen-Feldman \cite{Chen-Feldman15}),
shock diffraction in \S 3 (also see Chen-Xiang \cite{Chen-Xiang-1,Chen-Xiang-2}),
and the Prandtl-Meyer reflection described in \S 4 (also see Bae-Chen-Feldman \cite{BCF-13}).
In order to prove the convexity, we employ global properties of
admissible solutions, including the existence of the cone of monotonicity
discussed above.

\section{The Shock Reflection/Diffraction Problems
 and Free Boundary Problems for the Full Euler Equations}

When the vortex sheets and the deviation of vorticity become significant,
the full Euler equations are required.
In this section, we present mathematical formulation
of the shock reflection/diffraction problems for the full Euler
equations and the role of the potential theory for the shock
problems even in the realm of the full Euler equations.
In particular, the Euler equations for potential
flow, \eqref{1-r}--\eqref{1-p1},
are actually {\it exact} in an important region of the solutions
to the full Euler equations.

The full Euler equations for compressible fluids in
$\R^{3}_+:=\R_+\times \R^2, t\in\R_+:=(0,\infty), \x\in\R^2$, are of
the following form:
\begin{equation}\label{E-1}
\left\{\begin{aligned} &\partial_t\,\rho +\nabla_\x\cdot (\rho\vv)=0,\\
   %\qquad\qquad\qquad\qquad\qquad\qquad\qquad\mbox{(Euler 1755-59)}\\
&\partial_t(\rho\vv)+\nabla_\x\cdot\left(\rho \vv\otimes\vv\right)
  +\nabla p=0,\\
  %\quad \mbox{(Euler 1755-59, Cauchy 1827-29)}\\
&\partial_t(\frac{1}{2}\rho |\vv|^2+\rho
e)+\nabla_\x\cdot\big((\frac{1}{2}\rho|\vv|^2+\rho e + p)\vv\big)=0,
%\\  \qquad\quad \,\,\mbox{(Kirchhoff 1868)}
\end{aligned}
\right.
\end{equation}
where $\rho$ is the density, $\vv=(u, v)$ the fluid velocity, $p$
the pressure, and $e$ the internal energy. Two other important
thermodynamic variables are  the temperature $\theta$ and the energy
$S$.  The notation $\aaa\otimes\bb$ denotes the tensor product of
the vectors $\aaa$ and $\bb$.

Choosing $(\rho, S)$ as the independent thermodynamical variables,
then the constitutive relations can be written as
$(e,p,\theta)=(e(\rho,S), p(\rho,S),\theta(\rho,S))$ governed by
$$
\theta dS=de +pd\tau=de-\frac{p}{\rho^2}d\rho.
$$

For a polytropic gas, \begin{equation}\label{gas-1} p=(\gamma-1)\rho
e, \qquad e=c_v\theta, \qquad \gamma=1+\frac{R}{c_v},
\end{equation}
or equivalently,
\begin{equation}\label{gas-2}
p=p(\rho,S)=\kappa\rho^\gamma e^{S/c_v}, \qquad e=e(\rho,
S)=\frac{\kappa}{\gamma-1}\rho^{\gamma-1}e^{S/c_v},
\end{equation}
where $R>0$ may be taken to be the universal gas
    constant divided by the effective molecular weight of the particular
    gas,
    $c_v>0$ is the specific heat at constant volume,
    $\gamma>1$ is the adiabatic exponent, and $\kappa>0$ is any constant under
    scaling.

Notice that the corresponding three lateral Riemann problems in \S 2--\S 4
for system \eqref{E-1}
are all invariant under the self-similar scaling:
$(t, \x)\longrightarrow (\alpha t, \alpha \x)$ for any $\alpha\ne 0$.
Therefore, we seek self-similar solutions:
$$
(\vv, p,\rho)(t,\x)=(\vv,p,\rho)(\xi,\eta), \qquad
(\xi,\eta)=\frac{\x}{t}.
$$
Then the self-similar solutions are governed by the following
system:
\begin{equation}\label{Euler-sms-1}
\left\{\begin{aligned}
&(\rho U)_\xi +(\rho V)_\eta +2\rho=0,\\
&(\rho U^2+p)_\xi +(\rho UV)_{\eta}+3\rho
U=0,\\
&(\rho UV)_\xi +(\rho V^2+p)_{\eta}+3\rho
V=0,\\
&\big(U(\frac{1}{2}\rho q^2+\frac{\gamma p}{\gamma-1})\big)_\xi
+\big(V(\frac{1}{2}\rho q^2+\frac{\gamma p}{\gamma-1})\big)_\eta
+2(\frac{1}{2}\rho q^2+\frac{\gamma p}{\gamma-1})=0,
\end{aligned}
\right.
\end{equation}
where $q=\sqrt{U^2+V^2}$, and $(U,V)=(u-\xi, v-\eta)$ is the
pseudo-velocity.

The eigenvalues of system \eqref{Euler-sms-1} are
$$
\lambda_0=\frac{V}{U}\,\, \mbox{(repeated)},\qquad\,
\lambda_\pm=\frac{UV\pm c\sqrt{q^2-c^2}}{U^2-c^2},
$$
where $c=\sqrt{\gamma p/\rho}$ is the sonic speed.

When the flow is pseudo-subsonic, {\it i.e.},  $q<c$,
the eigenvalues $\lambda_\pm$ become complex and
thus the system consists of two transport equations
and two nonlinear equations of
elliptic-hyperbolic mixed type. Therefore, system
\eqref{Euler-sms-1}
is {\em hyperbolic-elliptic composite-mixed} in general.

The three lateral Riemann problems can be formulated
as the corresponding boundary value problems in the unbounded
domains. Then the boundary value problems can be further
formulated as the three corresponding free boundary problems.
The free boundary conditions are again the Rankine-Hugoniot
conditions along the free boundary $S$:
\begin{equation}\label{RH-full}
[L]_S=0, \quad
[\rho N]_S=0, \quad
[p+\rho N^2]_S=0,\quad
\big[\frac{\gamma p}{(\gamma-1)\rho}+\frac{N^2}{2}\big]_S=0,
\end{equation}
where $L$ and $N$ are the tangential and normal components of velocity $(U,V)$ along the free boundary,
that is, $|(U,V)|^2=L^2+N^2$.
The conditions along the sonic circles are the Dirichlet conditions for $(U,V, p,\rho)$ to be
continuous across the respective sonic circles.

\medskip
We now discuss the role of the potential flow equation \eqref{2-1}
in these free boundary problems whose boundaries also include
the fixed degenerate
sonic circles for the full Euler equations \eqref{Euler-sms-1}.

Under the Hodge-Helmoltz decomposition $(U,V)=D\varphi +W$ with
${\rm div} W=0$, the Euler equations \eqref{Euler-sms-1} become
\begin{eqnarray}
&&{\rm div}(\rho D\varphi)+2\rho=-{\rm div}(\rho W),\label{6.1}\\
&&D(\frac{1}{2}|D\varphi|^2+\varphi)+\frac{1}{\rho}Dp
=(D\varphi +W)\cdot DW+ (D^2\varphi+I)W,\label{6.2}\\
&&(D\varphi+W)\cdot D\omega +(1+\Delta \varphi)\omega=0,\label{6.3}\\
&&(D\varphi+W)\cdot DS=0,\label{6.4}
\end{eqnarray}
where $\omega={\rm curl}\, W={\rm curl} (U,V)$ is the vorticity of
the fluid, and $S=c_v\, {\rm ln} (p\rho^{-\gamma})$ is the entropy.

When $\omega=0, S=const.$ and $W=0$ on a curve $\Gamma$ transverse
to the fluid direction, we first conclude from \eqref{6.3} that, in
domain $\Omega_E$ determined by the fluid trajectories past $\Gamma$:
$$
\frac{d}{dt}(\xi,\eta)=(D\varphi+W)(\xi,\eta),
$$
we have
$$
\omega=0, \qquad \text{\it i.e.}\quad \text{curl}\, W=0.
$$
This implies that $W=const.$ since ${\rm div} W=0$. Then we
conclude
$$
W=0   \qquad \text{in}\,\, \Omega_E,
$$
since $W|_\Gamma=0$, which yields that the right-hand side of
equation \eqref{6.2} vanishes. Furthermore, from \eqref{6.4},
$$
S=const. \qquad\text{in}\,\, \Omega_E,
$$
which implies that
$$
p=const.\, \rho^\gamma.
$$
By scaling, we finally conclude that the solution of system
\eqref{6.1}--\eqref{6.4} in domain $\Omega_E$ is determined by
the Euler equations \eqref{1-r}--\eqref{1-p1}
for self-similar potential flow,
or the potential flow equation \eqref{2-1} with
\eqref{1-o}
for self-similar solutions.

For our problems in \S 2--\S 4,  we note that, in the supersonic states
joint with the sonic circles ({\it e.g.}
state (2) for {\bf Problem 2.3},
state (1) for {\bf Problem 3.3},
states $(0)$ and $(1)$ for {\bf Problem 4.3}),
\begin{equation}\label{across-sonic}
\omega=0, \quad W=0, \quad S=S_2.
\end{equation}
Then,  if our solution $(U,V,p,\rho)$ is $C^{0,1}$  and the gradient
of the tangential component of the velocity is continuous across the
sonic arc,
we still have \eqref{across-sonic} along $\Gamma_{\rm sonic}$ on the
side of $\Omega$. Thus, we have

\begin{theorem}\label{p-dominate}
Let $(U,V,p,\rho)$ be a solution of one of our problems, {\bf Problems 2.3, 3.3}, and {\bf 4.3},
such that
$(U,V,p,\rho)$  is $C^{0,1}$ in the open region formed by the reflected shock and the wedge boundary,
and the gradient of the tangential component of $(U,V)$ is continuous
across any sonic arc.
Let $\Omega_E$ be the subregion of $\Omega$ formed by the fluid
trajectories past the sonic arc, then, in
$\Omega_E$, the potential flow equation \eqref{2-1} with
\eqref{1-o} coincides with the full Euler equations
\eqref{6.1}--\eqref{6.4}, that is, equation \eqref{2-1} with
\eqref{1-o} is exact in the domain $\Omega_E$ for {\bf Problems
2.3, 3.3,} and {\bf 4.3}.
\end{theorem}

\begin{remark}
The regions such as $\Omega_E$ also exist in various Mach
reflection-diffraction configurations. {\rm Theorem \ref{p-dominate}}
applies to such regions whenever the solution $(U,V,p,\rho)$ is
$C^{0,1}$ and the gradient of the tangential component of $(U,V)$ is
continuous. In fact, {\rm Theorem \ref{p-dominate}}
indicates that, for the solution
$\varphi$ of \eqref{2-1} with \eqref{1-o}, the $C^{1,1}$--regularity
of $\varphi$ and the continuity of the tangential
component of the velocity field $(U,V)=\nabla\varphi$ are optimal
 across the sonic arc $\Gamma_{sonic}$.
\end{remark}

\begin{remark}
The importance of the potential flow equation \eqref{1-a} with
\eqref{1-b} in the time-dependent Euler flows even
through weak discontinuities
was also observed by Hadamard {\rm \cite{Hadamard}} through a different argument.
Moreover, for the solutions containing
a weak shock, the potential flow equation \eqref{1-a}--\eqref{1-b}
and the full Euler flow model \eqref{E-1} match
each other well up to the third order of the shock strength.
Also see
Bers {\rm \cite{Ber}},
Glimm-Majda {\rm \cite{Glimm-Majda}},
and Morawetz {\rm \cite{Morawetz}}.
\end{remark}

\section{Conclusion}

As we have discussed above, the three longstanding, fundamental transonic
flow problems can be all formulated as free boundary problems.
The understanding of these transonic flow problems requires our mathematical solution of these free boundary problems.
Similar free boundary problems also arise in many other transonic flow problems, including
steady transonic flow problems
including transonic nozzle flow problems (\textit{cf.} \cite{BF-11,Chen,Chen-Feldman07,LXY}),
steady transonic flows past obstacles (\textit{cf.} \cite{Chen,Chen-Chen-Feldman,C-Fang,CZZ,ChenS-Fang}),
supersonic bubbles in subsonic flow  (\textit{cf.} \cite{Cole-Cook,Morawetz-2}),
local stability of Mach configurations (\textit{cf.} \cite{ChenS-08}),
as well as
higher dimensional version of
{\bf Problem 2.3} (shock reflection-diffraction by a solid cone)
and {\bf Problem 4.3} (supersonic flow impinging onto a solid cone).
In \S 2--\S 5, we have discussed recently developed mathematical ideas, approaches, and techniques for solving these
free boundary problems.
On the other hand, many free boundary problems arising from transonic flow problems are still open and demand
further developments of new mathematical ideas, approaches, and techniques.

\bigskip
\medskip
{\bf Acknowledgements.} 
The authors would like to thank the Isaac Newton Institute for Mathematical Sciences,
Cambridge, for support and hospitality during the 2014 Programme on
{\it Free Boundary Problems and Related Topics} where work on this paper was undertaken.
They also thank Myoungjean Bae and Wei Xiang
for their direct and indirect contributions
in this paper. 
The work of Gui-Qiang G. Chen was supported in part by NSF Grant DMS-0807551,
the UK EPSRC Science and Innovation award to the Oxford Centre for Nonlinear PDE (EP/E035027/1),
the UK EPSRC Award to the EPSRC Centre for Doctoral Training
in PDEs (EP/L015811/1), and the Royal Society--Wolfson Research Merit Award (UK).
The work of Mikhail Feldman was
supported in part by the National Science Foundation under Grants
DMS-1101260, DMS-1401490, and by the Simons Foundation under the Simons Fellows Program.

\end{document}